\documentclass[a4paper,12pt]{article}
\usepackage{amsmath,amsfonts,amssymb,amscd}

\def\N{{\mathbb N}}
\def\Z{{\mathbb Z}}

\def\cH{{\mathcal H}}

\def\mfS{\mathfrak{S}}
\def\Y{\mathfrak{Y}}

\def\a{\alpha}
\def\b{\beta}
\def\g{\gamma}
\def\d{\delta}

\def\l{\lambda}

\def\moins{\raise 1pt\hbox{{$\scriptstyle -$}}}
\def\plus{\raise 1pt\hbox{{$\scriptstyle +$}} }
\def\s{\scriptstyle}
\def\bmoins{\fbox{\raise 1pt\hbox{{$\scriptstyle -$}}1}}

\def\Pfaff{\mathfrak{ P\hspace{-0.07 em}f\hspace{-0.05 em}a\hspace{-0.05 em}f%
\hspace{-0.07 em}f}}
\def\Pf{\mathfrak{ P\hspace{-0.07 em}f}}  
\def\Tab{\mathfrak{ T\hspace{-0.07 em}a\hspace{-0.05 em}b}}

\newtheorem{theorem}{Theorem}
\newtheorem{proposition}[theorem]{Proposition}
\newtheorem{lemma}[theorem]{Lemma}
\newtheorem{corollary}[theorem]{Corollary}
\def\proof{\noindent{\it Proof.}\ }

\def\QED{\hfill Q.E.D.}

\def\CARRE#1{\hbox{\vrule width \thickness
   \vbox to \carresize{\hrule height \thickness\vss
      \hbox to \carresize{\hss#1\hss}
   \vss\hrule height\thickness}
\unskip\vrule width \thickness}
\kern-\thickness}
\def\vsquare#1{\vbox{\CARRE{$#1$}}\kern-\thickness}

\def\smallyoung#1{%
  \newdimen\carresize \carresize=10pt%
  \newdimen\thickness \thickness=0.5pt%
  \vcenter{%
    \vbox{\smallskip\offinterlineskip%
      \halign{&\vsquare{##}\cr #1}}}}

\def\young#1{%
  \newdimen\carresize \carresize=16pt%
  \newdimen\thickness \thickness=0.5pt%
  \vcenter{%
    \vbox{\smallskip\offinterlineskip%
      \halign{&\vsquare{##}\cr #1}}}}

\newdimen\unit
\def\o{$\scriptscriptstyle{{\rm o}}$}
\def\grape(#1,#2)#3{\raise#2\unit\rlap{\kern#1\unit #3}\ignorespaces}

\def\gg{{\unit=1mm
\hbox {\grape(2,1.2){'}
       \grape(1,2)\o
       \grape(2,2)\o
       \grape(3,2)\o
       \grape(1.5,1)\o
       \grape(2.5,1)\o
       \grape(2,0)\o\
}\kern 3.5 \unit}}

\def\gfill{\leaders\hbox to 1.2em{\hss\gg\hss}\hfill}

\def\frise{\medskip \centerline{\hbox to 8cm{\gfill} }\bigskip}

\date{}
\begin{document}

\begin{center}
{\bf \large Pfaffians and Representations of the Symmetric Group}
\end{center}

\begin{center}
{\large Alain Lascoux}\\
 CNRS, IGM Universit\'e  de
  Marne-la-Vall\'ee\\
  77454 Marne-la-Vall\'ee Cedex, France\\
  Email: Alain.Lascoux@univ-mlv.fr\\
\end{center}

\frise 

\begin{abstract}
Pfaffians of matrices with entries $z[i,j]/(x_i+x_j)$, or determinants
of matrices with entries $z[i,j]/(x_i-x_j)$, where the antisymmetrical
indeterminates $z[i,j]$ satisfy the Pl\"ucker relations,
can be identified with a trace in an irreducible representation
of a product of two symmetric groups.
Using Young's orthogonal bases, one can 
write explicit expressions of such Pfaffians
and determinants, and recover in particular the  evaluation of Pfaffians
which appeared in the recent literature.

\end{abstract}

\frise

\vspace{0.2cm}



\vspace{0.2cm}

\section{Introduction}

Determinants or Pfaffians of order $n$ can be written in terms of 
the symmetric group $\mfS_n$. Determinants
can be considered as  generators of a 1-dimensional 
alternating representations.
But in the case of a determinant or a Pfaffian  
$$ \left| \frac{a_i-a_{j+n}}{x_i-x_{j+n}}  \right|_{1\leq i,j\leq n}
 \qquad ,\qquad 
\Pfaff\left(\, \frac{a_i-a_j}{x_i+x_j}\, \right)_{1\leq i<j\leq n} \, ,  $$
three symmetric groups occur~: $\mfS_n^a$ acts on the indeterminates
$a_i$, $\mfS_n^x$ acts on the $x_j$, and $\mfS_n^{ax}$ acts on the
indices of all indeterminates simultaneoulsy. 
This ``diagonal action'' satisfy a Cauchy-type property,
each irreducible representation of $\mfS_n^a$ occuring in the
expansion of the determinant, or of the Pfaffian, being tensored 
with a representation of $\mfS_n^x$ of conjugate type.

When $n=2m$ is even, since the space generated by the orbit of
the polynomial 
$(a_1 \moins a_2)(a_3\moins a_4)\cdots (a_{n-1}\moins a_n)$ 
under $\mfS_n^a$ is a copy of the irreducible representation 
$V_{[m,m]}^a$ of type $[m,m]$,
 this forces the Pfaffian to lie in the space 
   $V_{[m,m]}^a \otimes V_{[2,\ldots,2]}^x$.

An easy analysis shows that moreover the Pfaffian 
$\Pfaff\left(\frac{a_i-a_j}{x_i+x_j}  \right)$ is diagonal in 
Young's orthogonal basis (and thus, can be considered as a trace).  
In fact, the same analysis remains
valid  (this is our main theorem, Th.\ref{th:Pfaffzg}) in the more general case 
$$ \Pfaff\left( z[i,j]\, g[i,j]  \right)_{1\leq i<j\leq n} \, ,  $$
when taking  antisymmetric indeterminates
$z[i,j]=\moins z[j,i]$ 
satisfying the Pl\"ucker relations 
(we say \emph{Pl\"ucker indeterminates}), instead of $(a_i-a_j)$,
and symmetric indeterminates $g[i,j]=g[j,i]$ instead of 
$(x_i+x_j)^{-1}$. 
For specific $z[i,j]$ and $g[i,j]$, one may be able to write another 
element belonging to the same representation. Checking that two elements
in the same irreducible representation coincide is very
easy, and reduces to compute some specializations. 

The most general case that we consider is 
$\Pfaff\left( a[i,j]b[i,j] z[i,j]^{-1} \right)$, with three families
of Pl\"ucker indeterminates. In that case the Pfaffian factorizes
in two factors separating the $a[i,j]$'s and $b[i,j]$'s 
(Th.\ref{th:Pfaffabz}).

A connection with the theory of symmetric functions is provided 
by specializing the Pl\"ucker indeterminates into
$S_\l(A+x_i +x_j)(x_i-x_j)$, $S_\l(A)$ being a fixed Schur function,
to the alphabet of which one adds the letters $x_i,x_j$ (\cite{Cbms}).
Thus Th.\ref{th:Pfaffabz} gives the factorization of 
$$ \Pfaff\left( \frac{ S_\l(A+a_i+a_j) 
S_\mu(B+b_i+b_j)}{S_\nu(Z+z_i+z_j)}
 \frac{(a_i\moins a_j)(b_i\moins b_j)}{(z_i-z_j)}          \right)  \, ,$$
for three Schur functions, and three families of indeterminates.

In the case of $\Pfaff\left(\frac{a_i-a_j}{x_i+x_j}  \right)$ first
considered by Sundquist \cite{Sundquist}, 
and that we have taken as our generic case, 
it is easy to write a
determinant which also lies in the space 
$V_{[m,m]}^a \otimes V_{[2,\ldots,2]}^x$. 
Specializing half of the $a_i$'s to $1$, the others to $0$, 
one recovers the determinantal expression of Sundquist for this
Pfaffian (Th.\ref{th:Sundquist}).

Ishikawa \cite{Ishikawa}, Okada \cite{Okada}, 
M. Ishikawa, S. Okada, H. Tagawa and J. Zeng
\cite{IshiOkada} have given different 
generalizations of Sundquist's Pfaffian. We show 
how to connect  their results to Th.\ref{th:Pfaffzg} and 
Th.\ref{th:Pfaffabz}.

In section 8, we go back to determinants, and show 
how to relate 
\begin{equation*}
\left| \frac{z[i,j]}{x_i^2-x_j^2}    \right|_{1\leq i\leq m<j\leq n}
\quad \text{and}\quad \Pfaff\left(\frac{z[i,j]}{x_i+x_j}  \right) \, ,
\end{equation*}
the indeterminates $z[i,j]$ still satisfying the Pl\"ucker relations.
A corollary of this analysis is that the above determinant
is, up to straightforward factor, symmetrical in $x_1,\ldots, x_n$,
and not only symmetrical in $x_1,\ldots, x_m$ and
$x_{m+1},\ldots, x_{2m}$ separately. 
Some determinants \quad  $\det\Bigl( S_\l(A+x_i+x_j)\\
 S_\mu(B+x_i+x_j)^{-1}  \Bigr)$
present a special interest in the theory of orthogonal polynomials,
or of the six-vertex model.

To be self-contained, and for lack of a reference appropriate to our
needs, we first recall some properties of representations.  
In the last section, we give more details about the polynomial
bases that one deduces from Young's orthogonal idempotents.

\section{Representations of the symmetric group}

\subsection{Young's idempotents}

The group algebra $\cH$
of the symmetric group $\mfS_n$ has by definition 
a linear basis consisting of all the permutations of $1,2,\ldots,n$.

Young described another basis $e_{tu}$, indexed by pairs of
standard tableaux of the same shape with $n$ boxes. 
These elements are 
\emph{matrix units}, in the sense that they satisfy the relations
\begin{eqnarray}
  e_{t,u} e_{u,v} &=& e_{t,v}  \, ,  \\
  e_{t,u} e_{w,v} &=& 0 \quad \text{if}\ w\neq u\, . 
\end{eqnarray}

In particular, the $e_{t,t}$ are idempotents: $e_{t,t} e_{t,t}= e_{t,t}$,
and the identity decomposes as
$$ 1 = \sum_t e_{t,t} \, , $$
where the sum is over all standard tableaux of $n$ boxes.  
The subsum 
\begin{equation}   \label{IdempCentral}
e_\l=\sum_{t\in Tab(\l)}   e_{t,t} 
\end{equation} 
over  standard tableaux of a given
shape $\l $ is the \emph{central idempotent} 
of index $\lambda$.

\subsection{Specht representations}

Given any $t$, the right module $e_{t,t}\, \cH$ is an irreducible
representation of the symmetric group, 
with basis $\{ e_{t,u}:\, u $ has the same shape as $t\}$. 

There are simpler models of irreducible representations, in particular
spaces of polynomials which are called \emph{Specht representations}, 
though they have been defined by Young\footnote{
Young \cite[Theorem IV, p.591]{Young} uses the  picturesque terminology
"has the same substitutional qualities", to state that 
the space $e_{t,t}\, \cH$ is isomorphic to the space generated
by some products of Vandermonde determinants.}.

Bases are still indexed by standard tableaux of a given shape,
but now tableaux are interpreted as polynomials as follows.

A column tableau $\smallyoung{\s k\cr\s j\cr\s i\cr }$ is interpreted as
the Vandermonde determinant in  the variables $x_i,x_j,\ldots, x_k$~:
$$  \young{k\cr j\cr i\cr} = (x_i-x_j)\, (x_i-x_k)\, (x_j-x_k) 
:= \Delta^x(i,j,k) \, ,$$
and a tableau stands for the product of its columns~:
$$\young{ 5\cr 3 &6\cr 1&2&4\cr}  = 
\Delta^x(1,3,5)\, \Delta^x(2,6)  \, \Delta^x(4)\, .    $$

We shall denote this polynomial $\Delta_t^x$ and call it 
\emph{Specht polynomial}. 
The orbit of any $\Delta_t^x$  under the symmetric group
(permuting the variables $x_i$) has $n!$ elements, whose linear 
span is of dimension the number of standard tableaux of the same
shape as $t$.

More precisely, Young obtained, in the case of zero characteristic~:

\begin{proposition}  \label{th:Specht1} 
 Given a partition $\l$, the linear span of the 
polynomials $\Delta_t^x$, $t$ varying over the set $\Tab(\lambda)$ 
of standard tableaux of shape $\l$,
is an irreducible representation of the symmetric group.
\end{proposition}

Fixing a shape $\lambda$, there are two ``extreme tableaux'':
the one such that its columns are filled with consecutive letters,
and that we shall call \emph{top tableau} and denote $\zeta$. 
For shape $[2,3,4]$, the top tableau is 
$$ \zeta = \young{3 &6\cr 2 & 5 & 8\cr 1&4 &7 &9\cr}  $$
and gives the Specht polynomial
$$ \Delta_\zeta^x = 
  \Delta^x(1,2,3)\, \Delta^x(4,5,6)\, \Delta^x(7,8)\, \Delta^x(9)\, .  $$ 

Similarly, the \emph{bottom tableau} has its rows filled with
consecutive letters. We denote it by $\aleph$: 
$$ \aleph = \young{8 &9\cr 5&6&7\cr 1&2&3&4 \cr} \, $$
with Specht polynomial
$$  \Delta_\aleph^x = 
 \Delta^x(1,5,8)\, \Delta^x(2,6,9)\, \Delta^x(3,7)\, \Delta^x(4)\, .  $$

The standard tableaux of a given shape may be generated by using
simple transpositions, starting with $\zeta$. By the notation 
$\Tab(\l)$ we mean this ranked poset, with top element $\zeta$
and bottom one, $\aleph$.  
The \emph{distance} $\ell(t,u)$ of two tableaux is the distance 
in  $\Tab(\l)$.

The decomposition of any element of the Specht representation
in the basis $\Delta_t^x$ is given by a so-called \emph{straightening 
algorithm} (cf. \cite{DKR, Desarmenien, Turbo}). 

We shall need only one coefficient in such an expansion,
the coefficient of $\Delta_\aleph^x$.
Given a tableau $t$ with $n$ boxes, and a function $f(x_1,\ldots, x_n)$, 
let us write $f(t)$ for the specialization where each $x_i$ is specialized
to $r$ if $i$ lies on row $r$ (rows are numbered from the bottom,
starting with $0$).

\begin{lemma}   \label{th:Specht2}  
Given a partition $\lambda$, and a linear combination 
$f(x_1,\ldots, x_n)= \sum_t c_t\, \Delta_t^x $, 
with coefficients $c_t$ independent of $x_1,\ldots, x_n$, then the 
coefficient $c_\aleph$ is equal to 
$$   f(\aleph)/ \Delta_\aleph^x(\aleph)  \ .$$
\end{lemma}

\proof All other tableaux than $\aleph$  
 have in some column two entries which lie in the same row 
of $\aleph$.                          \QED     

The Specht representation can occur in many disguises.
Let us call \emph{Pl\"ucker indeterminates} 
anti-symmetric indeterminates $z[i,j]=-z[j,i]$,
satisfying Pl\"ucker relations\cite{BL} for all quadruples
of different integers~:
$$  z[i,j] z[k,l] -z[i,k]z[j,l] + z[j,k]z[i,l]= 0 \ . $$

A  typical example is obtained by taking a $2\times \infty$ generic
matrix $M$, and defining $z[i,j]$ to be the minor on 
columns $i,j$ of $M$. More generally, one takes an $N\times\infty$
generic matrix $M$, one chooses $N-2$ columns of index $\a,\b,\ldots$
and define $z[i,j]$ to be the minor of maximal order of $M$
on columns $i,j,\a,\b,\ldots$, with $i,j\neq \a,\b,\ldots$.

The following proposition gives another description of Specht 
representations for shape $[m,m]$.

\begin{proposition}  \label{th:Specht3}
 Given an even positive number $n=2m$,
 let $z[i,j]$, $1\leq i, j \leq n$ be Pl\"ucker indeterminates.
Given any numbering $u$ of the boxes of the diagram $[m,m]$,
let $z[u]$ be the product of all $z[i,j]$, where $[i,j]$ is a column
of $u$.Let the symmetric group $\mfS_n$ act on the variables $z[i,j]$ by
permutation of $1,2,\ldots,n$.

Then the correspondence $z[u]=\prod z[i,j] \to \prod (x_i-x_j) $
induces an isomorphism of representations of $\mfS_n$.
In particular, $\{ z[t]\}$, where $t$ runs over all standard
tableaux of shape $[m,m]$, is a linear basis of the
span of all $z[u]$.
\end{proposition}

In short, when one has 
$$  z[1,2]z[3,4] -z[1,3]z[2,4] +z[1,4]z[2,3] =0 \, $$
one can as well read 
$$  (a_1-a_2)(a_3-a_4) - (a_1-a_3)(a_2-a_4) + (a_1-a_4)(a_2-a_3) =0 $$
or
$$  \young{ 2 &4\cr 1&3\cr} -  \young{ 3 &4\cr 1&2\cr}
 +  \young{ 4 &3\cr 1&2\cr} = 0  $$
without loss of generality.

We shall use, for a rectangular shape with two columns
of length $m=n/2$, three models of representations. 
The first one is the usual Specht representation,
 generated by the action of $\mfS_n$ on 
$$ \Delta^x(1,2,\ldots,m)\,  \Delta^x(m\plus1,\ldots,n)\, .   $$

The second model is the image of the first one under the correspondence
$(x_i-x_j) \to z[i,j]$. The Specht polynomial corresponding to
the top tableau of shape $2^m$ will now be 
$$  \Y^{z[\,]}(\zeta) := 
\Delta^{z[\,]}(1,2,\ldots,m)\,  \Delta^{z[\,]}(m\plus1,\ldots,n) 
= \prod_{1\leq i<j\leq m} z[i,j]\, \prod_{m+1\leq i<j\leq n} z[i,j]
\, .  $$

For the third one, one starts with a symmetric matrix 
$G$, with entries $g[i,j]=g[j,i]$.
Let us denote the minor consisting of rows $i_1,\ldots, i_m$ and
columns $i_{m+1},\ldots, i_n$ by 
$$   g[ i_1,\ldots, i_m\, |\, i_{m+1},\ldots, i_n]  \, .$$
The symmetric group $\mfS_n$ acts formally by permuting the indices 
of such minors.
Now, Kronecker \cite{Kronecker, Muir}  
has shown that such minors satisfy the Pl\"ucker relations\footnote{
and, of course, all the relations obtained by permuting
the rows and the columns of the original matrix,
in such a way as to obtain another symmetric matrix.
We could reprove directly Kronecker's relations
by introducing extra variables $a_1,a_2,\ldots$
and evaluating the Pfaffian of the antisymmetric matrix 
$\bigl[(a_i-a_j)g[i,j] \bigr]_{i,j=1,\ldots,n}$,
as will become clear later.}
$$ \sum_{i=0}^m   (-)^i\, 
 g[1,\ldots,m\moins1, m\plus i\, |\, m\plus1,\ldots, 
  \widehat{m\plus i},\ldots,  n] =0 \, . $$

This implies the following proposition~:
\begin{proposition}
Let $G$ be a symmetric matrix of order $n=2m$. 
Then the linear span of the minors in the orbit of
$$  g[1,\ldots,m\, |\, m\plus1,\ldots,n] $$
under permutation of indices, is an irreducible
representation of $\mfS_n$ of shape $2^m$.
\end{proposition}

Using the correspondence 
$$ g[ i_1,\ldots, i_m\, |\, i_{m+1},\ldots, i_n]
\leftrightarrow \Delta^x(i_1,\ldots, i_m)\, \Delta^x(i_{m+1},\ldots, i_n) $$
one can still speak of a \emph{Specht basis} for these  minors
of a symmetric matrix.

For example, for $m=3$, the space has basis
$$ g[123\, |\, 456],\, g[124\, |\, 356],\,
g[125\, |\, 346],\,  g[134\, |\, 256],\,
g[135\, |\, 246],\, ,$$
and one can directly check the relation
$$  g[123\, |\, 456] -g[124\, |\, 356] +g[125\, |\, 346] 
- g[126\, |\, 345]  =0 \, .$$

In detail, for $m=2$,
the sum $ g[12|34]-g[13|24]+g[14|23] $ expands into 
$$ 
g[1, 3] g[2, 4] - g[1, 4] g[2, 3] - g[1, 2] g[3, 4] + g[1, 4] g[3, 2]
     + g[1, 2] g[4, 3] - g[1, 3] g[4, 2]
$$
which is indeed zero, because $g[i,j]=g[j,i]$.

Pl\"ucker relations, hence Specht representations, also occur
in the theory of symmetric functions.

Indeed, given two alphabets\cite{Cbms} $A=\{a\}$, $B=\{b\}$, the 
\emph{complete functions} $S_k(A-B)$ are defined by the generating
function
$$ \sum_k z^k S_k(A-B) = \prod_{b\in B} (1-zb) 
    \prod_{a\in A} (1-za)^{-1} \, ,$$
putting $S_{k}=0$ for $k<0$.  The \emph{Schur function}
$S_v(A-B)$, $v\in \Z^r$, has the determinantal expression    
$\det\left( S_{v_j+j-i}(A-B) \right)$. 
When $B$ is the two-letters alphabet $B=\{x,y\}$, then 
$(x-y)S_v(A-B)$, denoted $(x\moins y) S_v(A -x-y)$, 
is equal to the maximal minor 
on columns $v_1,v_2\plus 1, \ldots, v_r \plus r\moins 1, x,y$ of   
$$  \begin{array}{r} \s columns \\ \phantom{S_0(A)}\\
\phantom{S_0(A)}\\ \phantom{\vdots}\\  \phantom{S_0(A)}  \end{array}
  \begin{bmatrix} 
   \s 0 & \s 1 &\s 2 & \s 3 & \cdots &   & \s x & \s y \\
  \noalign{\kern 3pt \hrule \kern 3pt}
  S_0(A) & S_1(A) & S_2(A) &S_3(A) &\cdots & & x^{r+1} & y^{r+1} \\
  S_{-1}(A) & S_0(A) & S_1(A) & S_2(A)  &\cdots & & x^r & y^r \\
   \vdots  & \vdots  &\vdots  &\vdots  &        & &\vdots &\vdots \\
  S_{-r-1}(A) & S_{-r}(A) & S_{-r+1}(A) & S_{-r+2}(A)  &\cdots & & 1 & 1  
\end{bmatrix} \, .
$$

Therefore, for a given $v\in\Z^r$, and a given alphabet $A$,
the $z[i,j] = (x\moins y) S_v(A -x-y)$ (resp. 
 $z[i,j] = (x\moins y) S_v(A +x+y)$) satisfy the Pl\"ucker relations.

\subsection{Orthogonal representations}   

Given a shape $\lambda$, then 
$ \Delta_\zeta^x\, e_{\zeta,\zeta} = \Delta_\zeta^x$, and 
the polynomials 
 $\Delta_\zeta^x\, e_{\zeta,t}$, for $t\in \Tab(\l)$, 
constitute another basis of the Specht representation. 
Taking a linear order compatible with the poset structure, then
the matrix of change of basis is lower triangular.
A more precise information is given by using the Yang-Baxter
relations (see the last section). 

Let us call \emph{Young's basis} the basis
proportional to $\Delta_\zeta^x\, e_{\zeta,t}$, such that the leading
term (with respect to the poset $\Tab(\l)$) of each 
$ \Y(t)$ be $\Delta_t^x$, and  denote it
$\{ \Y(t) :\, t\in \Tab(\l)\}$.

For example, for shape $[3,3]$, the poset of tableaux is
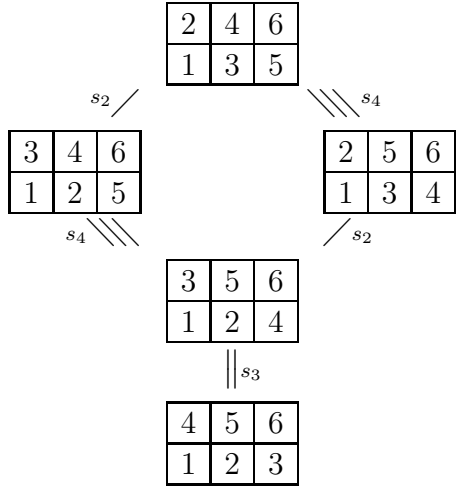
\begin{equation*}   
\begin{array}{rcl}   
          & \young{2&4&6\cr 1&3&5\cr}                \\[4pt]
  \raise 3pt \hbox{$\s s_2$}  \diagup &
            & \hskip -6pt\diagdown\hskip -6pt\diagdown
         \hskip -6pt\diagdown \raise 3pt \hbox{$\s s_4$} \\{}
\young{3&4&6\cr 1&2&5\cr} & &  \young{2&5&6\cr 1&3&4\cr}\\[10pt]
{\s s_4} \diagdown\hskip -6pt\diagdown\hskip -6pt\diagdown &&
   \diagup {\s s_2} \\
               & \young{3&5&6\cr 1&2&4\cr} &  \\[13pt]
           & \phantom{\s s_3} \big\| {\s s_3} &\\[2pt]
                       & \young{4&5&6\cr 1&2&3\cr}
 \end{array}
\end{equation*}
and the matrix expressing Young's basis in terms of
the Specht basis (reading successive rows) is
$$\left [\begin {array}{ccccc} 1&0&0&0&0\\\noalign{\medskip}-1/2&1&0&0
&0\\\noalign{\medskip}-1/2&0&1&0&0\\\noalign{\medskip}1/4&-1/2&-1/2&
1&0\\\noalign{\medskip}2/3&-1/3&-1/3&-1/3&1\end {array}\right ]\,  
=\, 
\left [\begin {array}{ccccc} 1&0&0&0&0\\\noalign{\medskip}1/2&1&0&0&0
\\\noalign{\medskip}1/2&0&1&0&0\\\noalign{\medskip}1/4&1/2&1/2&1&0
\\\noalign{\medskip}-1/4&1/2&1/2&1/3&1\end {array}\right ]^{-1} 
\, .  $$

In that special case, the two allowable linear orders on the graph
give the same matrices, we did not need numbering the tableaux.

To handle other models of irreducible representations, 
we first need to characterize the elements corresponding to
$\Delta_\zeta^x$ in these models.  In fact, 
$\Delta_\zeta^x$ is the only polynomial in the Specht representation,
such that $\left( \Delta_\zeta^x\right)^{s_i} =- \Delta_\zeta^x$
for all $i$ such that $i,i+1$ are in the same column of $\zeta$.
This property still characterizes a unique  element 
in a different copy of the Specht representation,
that we shall still denote  $ \Y(\zeta)$.

We shall still call \emph{Young basis} the basis proportional to
 $\{ \Y(\zeta)\, e_{\zeta,t}\}$,
with the same factors of proportionality as in the case where 
we start with $\Delta_\zeta^x$. 
More details about how
to compute such a basis are given in the last section.

\subsection{Cauchy Formula}

Let us now take  two symmetric groups $\mfS_n^x$, $\mfS_n^a$
acting respectively on the variables $x_1,\ldots, x_n$,
and $a_1,\ldots, a_n$, with generators $s_i^x$, $s_i^a$.

We also use the group $\mfS_n= \mfS_n^{ax}$,
 which permutes simultaneously
the variables $a_i$ and $x_i$.
We write its generators $s_i$, instead of $s_i^{ax}$. These are such that 
$$   s_i =  s_i^a s_i^x = s_i^x s_i^a \, .$$

We shall give a decomposition of the two 1-dimensional idempotents \\
$\square:= (n!)^{-1}\sum \sigma$ and 
$\nabla:=(n!)^{-1} \sum (\moins 1)^{\ell(\sigma)} \sigma$,
with respect to the groups $\mfS_n^x$ and  $\mfS_n^a$.

Indeed, by definition 
$$ \square^{ax} = \frac{1}{n!}\, \sum_\sigma \sigma^x \sigma^a\ , $$
with pairs of permutations in $\mfS_n^a$, $\mfS_n^x$
commuting with each other.

Taking the basis of Young idempotents, instead of the basis
of permutations,  one obtains  
a Cauchy-type formula~:
\begin{equation}  \label{CauchyCarre}
  \square^{ax} = \sum_{t,u}  \frac{1}{d(t)} e_{t,u}^x\, e_{t,u}^a    \ ,
\end{equation}
where the sum is over all pairs of standard tableaux of the same shape,
and where $d(t)$ is the number of tableaux of the same shape as $t$.

The element $\nabla^{ax}$ is obtained from $\square^{ax}$ under the 
isomorphism induced by~:
$$ s_i^x \to \widehat{s_i^x}:= - s_i^x \ , \ s_i^a \to  s_i^a   $$


From the expressions of $\square$ and $\nabla$ as sums of
products of permutations, one sees
that for any simple transposition, 
\begin{eqnarray}
 \square^{ax} \, s_i^x = \square^{ax} \, s_i^a \  &\& &
  s_i^x\, \square^{ax}  = s_i^a\, \square^{ax}  \\
 \nabla^{ax} \, s_i^x = -\nabla^{ax} \, s_i^a \  &\& &
  s_i^x\, \nabla^{ax}  = -s_i^a\, \nabla^{ax}  \ .
\end{eqnarray}

By taking products, one gets 
$$  \square^{ax} \, s_i^x s_j^x \cdots s_k^x =
 \square^{ax} \, s_i^a s_j^x \cdots s_k^x 
     = \square^{ax} \, s_j^x \cdots s_k^x s_i^a =\cdots \ , $$
and therefore, for any permutation $\sigma$, 
\begin{eqnarray}  \square^{ax} \, \sigma^x = \square^{ax} \, (\sigma^a)^{-1}  
 \  &\& & \ \sigma^x\, \square^{ax} = (\sigma^a)^{-1} \, \square^{ax}  \\
\nabla^{ax} \, \sigma^x = (\moins 1)^{\ell(\sigma)}
                       \nabla^{ax} \, (\sigma^a)^{-1} 
\ &\& &  \ \sigma^x\, \nabla^{ax} = (\moins 1)^{\ell(\sigma)}
  (\sigma^a)^{-1}\, \nabla^{ax} \ . 
\end{eqnarray}


Restricting  $\square^{ax}$ or $\nabla^{ax}$ 
 to a representation of $\mfS_n^a$ of type
  is achieved by multiplying 
$\square^{ax}$ or $\nabla^{ax}$ by $e_\l^a$. 

Since an idempotent $e_{t,t}^x$ is sent under $s_i^x\to -s_i^x$ to
$e_{t^\sim,t^\sim}^x$, where $t\to t^\sim$ denotes the transposition
of tableaux, using expression 
(\ref{IdempCentral}), we obtain the two equivalent expansions~:
\begin{eqnarray}
 \square^{ax}\, e_\l^a &=& \sum_{t,u\in \Tab(\l)} 
                     \frac{1}{d(t)}\,   e_{t,u}^a\, e_{t,u}^x  \ , \\
 \label{DecompoNabla}
\nabla\, e_\l^a &=& \sum_{t,u\in \Tab(\l)}
     (-1)^{\ell(t,u)} \frac{1}{d(t)}  e_{t,u}^a\, e_{t^\sim,u^\sim}^x  \ .
\end{eqnarray}

For example, there are $4$ standard tableaux with three boxes~:
$$ \a= \smallyoung{\s1& \s2&\s3\cr}\ ,\ 
   \b=\smallyoung{\s2\cr \s1&\s3\cr}\ ,\
   \g=\smallyoung{\s3\cr  \s1&\s2\cr}\ ,\
   \d=\smallyoung{\s3\cr \s 2\cr \s 1\cr} \ ,$$
and the elements $\square^{ax}$ and $\nabla^{ax}$ decompose as ~:
\begin{eqnarray*}
 \square^{ax} &=& e_{\a,\a}^a e_{\a,\a}^x + \frac{1}{2}\left(
   e_{\b,\b}^a e_{\b,\b}^x  + e_{\b,\g}^a e_{\b,\g}^x 
+ e_{\g,\b}^a e_{\g,\b}^x +e_{\g,\g}^a e_{\g,\g}^x  \right)
+ e_{\d,\d}^a e_{\d,\d}^x   \\  
\nabla^{ax}  &=& e_{\a,\a}^a e_{\d,\d}^x + \frac{1}{2}\left(
   e_{\b,\b}^a e_{\g,\g}^x  - e_{\b,\g}^a e_{\g,\b}^x
- e_{\g,\b}^a e_{\b,\g}^x +e_{\g,\g}^a e_{\b,\b}^x  \right)
+ e_{\d,\d}^a e_{\a,\a}^x\, ,
\end{eqnarray*}
the middle part being the component of type $[1,2]$.

\section{Pfaffians}

Let $Z=\bigl[ z[i,j] \bigr]$ be an anti-symmetric matrix of even order
$n$.  Its determinant is the square  of a  
function of the $z[i,j]$, which is called
the \emph{Pfaffian} of $Z$. 
The Pfaffian is a certain sum, with coefficients $\pm$,
of products $z[i,j]\cdots z[k,l]$.
We refer to \cite{Knuth} for an historical and complete presentation.

Deciding to write each monomial in $z[i,j]$ 
 according to some lexicographic order on the variables,
one can erase $z$, brackets and commas,  and use permutations~:
 $$ z[i,j]\, z[k,l]\, z[p,q] \to  [i\, j\,\, k\, l\,  p\, q]  $$

The Pfaffian of $Z$ has become an alternating sum of permutations
which can be defined recursively as follows \cite{Knuth}.
For any vector $v\in \N^n$ of even length $n$, let 
$$\Pf(v) = \sum_{i=1}^{n-1} (-1)^i\,  \Pf(v\setminus \{v_i,v_n\})\cdot
   [v_i,v_n]  \ , $$
where the product is the concatenation product,
and where $v\setminus \{v_i,v_n\}$ means suppressing the
 components $v_i,v_n$ inside $v$.

The initial case is $\Pf([\, ]) = [\, ]$. \\
$\Pf([1,2])=[1,2] \quad ; \quad 
 \Pf([1,2,3,4])=[1,2,3,4]-[1,3,2,4]+[2,3,1,4] \ ,$\\ 
$ \Pf([1,2,3,4,5,6])= 
[1,2,3,4,5,6]-[1,2,3,5,4,6] +[1,2,4,5,3,6]\\
\hspace*{10pt} -[1,3,2,4,5,6]+[1,3,2,5,4,6]-[1,3,4,5,2,6] -[1,4,2,5,3,6]\\
\hspace*{10pt} +[1,4,3,5,2,6]+[2,3,1,4,5,6] -[2,3,1,5,4,6]+[2,3,4,5,1,6]\\
\hspace*{10pt} +[2,4,1,5,3,6] -[2,4,3,5,1,6]-[3,4,1,5,2,6]+[3,4,2,5,1,6]\, . $

Let $\Pf_n := \Pf([1,\ldots,n])$ be the above sum of permutations 
(we use the notation $\Pfaff(Z)$ for the Pfaffian of an antisymmetric
matrix, and $\Pf_n$ for the formal sum of permutations).
Our data, the $z[i,j]$, were such that $z[i,j]=-z[j,i]$, and that
$z[i,j] z[k,l]=z[k,l]z[i,j]$.
We shall see in the next
proposition,whose proof is immediate,  that 
the permutations appearing in $\Pf_n$ are cosets representatives,
modulo the symmetries possessed by the $z[i,j]$.

Indeed, let $\mfS_{n/2}$ be the symmetric group which permutes
the blocks $[1,2],[3,4],\ldots, [n\moins 1,n]$, and
let $\Theta_n$ be the sum of its elements.
Simple transpositions $s_1,s_3,\ldots, s_{n-1}$ commute with
$\Theta_n$.


\begin{proposition}   \label{Pfaff2Sg}
For even $n$, the alternating sum of all permutations can be
factorized as follows~: 
\begin{eqnarray*}
n!\, \nabla &=& (1-s_1)(1-s_3)\cdots (1-s_{n-1})\, \Theta_n\, \Pf_n  \\
   &=& \Theta_n\, (1-s_1)1-s_3)\cdots (1-s_{n-1})\, \Pf_n\, . 
\end{eqnarray*}

\end{proposition}

As a consequence, we may write the Pfaffian of $Z$ as
$$ \Pfaff(Z) =   z[1,2]\, z[3,4]\cdots z[n\moins 1,n]\, \nabla 
 \, \frac{n!}{2^n (n/2)!} \ , $$
since $s_1,s_3,\dots$ and the permutations in $\Theta_n$ act trivially
on the product $z[1,2]z[3,4]\cdots$.

For example, for $n=6$, 
\begin{multline*}
  \sum _{\sigma\in\mfS_6} (\moins 1)^{\ell(\sigma)}\sigma =
      (1-s_1)(1-s_3)(1-s_5) \bigl([123456]+[125634]+[341256]   \\
+[345612]+[561234]+[563412] \bigr)\, 
\Pf_6\, .
\end{multline*}

An interesting approach, due to Luque and Thibon \cite{LuqueThibon},
 to combinatorial properties of Pfaffians
is through shuffle algebras. As a matter of fact, the same methods
give also the \emph{Hafnian}, i.e. the image of the Pfaffian (as an
element of the group algebra) under the involution $s_i\to -s_i$,
$i=1,\ldots, n\moins 1$. We shall not need this approach, having written
the Pfaffian in terms of the alternating sum of all permutations.  

In \cite{LLT} one finds Pfaffians and determinants associated
to a family of formal series, which are needed in geometry.

\section{Pfaffian, with the help of two symmetric groups}   

The main case that we want to treat now is the case of an 
antisymmetric matrix with entries 
$$ z[i,j]\, g[i,j] \, , $$
the $z[i,j]$ satisfying the Pl\"ucker relations, and the 
$g[i,j]$ being symmetrical: $g[i,j]=g[j,i]$.

Of course, any antisymmetric matrix $N= \Bigl[ n[i,j]\Bigr]$ 
can be written in this way,
introducing extra variables $a_i$, and writing
$$ n[i,j] = (a_i-a_j)\, \frac{n[i,j]}{a_i-a_j} \ .$$

The Pfaffian, being a sum of products of $z[i,j]$,
 belongs to the irreducible representation of shape 
$[m,m]$ 
(with respect to the symmetric group $\mfS_n^z$
acting on the indices of the indeterminates $z[i,j]$,
$i=1,\ldots,n$, $m=n/2$). 
Hence, it can be expressed as a linear combination
of Specht elements~:
$$ c_\zeta(g) z[\zeta] +\cdots +c_{\aleph}(g) z[\aleph] \ . $$

Thanks to Prop.\ref{th:Specht3}, the coefficients 
$c_t(g)$ are the same as for the 
specialization $z[i,j]= (a_i-a_j)$.

The case of the Pfaffian of $ (a_i-a_j)\, (x_i+x_j)^{-1}$ 
has been treated by Sundquist \cite{Sundquist},
but we need a more extensive description than his.
 
The coefficient $c_{\aleph}(g)$ is obtained by 
specializing $a_1=1 = \cdots = a_m$,
 $a_{m+1} =0=\cdots = a_n$. 
In that case, the sum over all permutations in $\mfS_n^{a,g}$
$$ \sum (-)^{\ell(\sigma)} \left(  (a_1\moins a_2) g[1,2]\, 
  (a_3\moins a_4) g[3,4]\,  \cdots \right)^\sigma  $$
reduces to 
$$ \sum (-)^{\ell(\sigma)} 
   \bigl(   g[1,m\plus 1]\, g[2,m\plus 2]\, \cdots 
    g[m\moins 1, n] \bigr)^\sigma  \ ,$$
where the sum is now only over the subgroup 
$ \mfS_{m} \times \mfS_{m}$ which permutes 
$1,\ldots, m$, and $m\plus 1,\ldots, 2m$ separately.  

This sum is equal to
\begin{equation}  \label{Detg}
   m!\, g[1\ldots m\, \,|\, m\plus 1\ldots n]  \, .
\end{equation}
Therefore the Pfaffian of the matrix $\Bigl[(a_i-a_j)g[i,j] \Bigr]$
is, up to a normalization constant, equal to the image of
$$ \Omega_{a,g}:= (a_1-a_{m+1})\cdots (a_m-a_n)\, 
                    g[1\ldots m\, \,|\, m\plus 1\ldots n]  $$
under the anti-symmetrization $\nabla^{ag}$.

This last element belongs to the Specht representation of $\mfS_n^a$ of shape
$[m,m]$, and therefore, thanks to (\ref{DecompoNabla}), belongs to the 
space\footnote{We do not need to know Kronecker's relations. 
Having only a component of type $[m,m]$ for $\mfS_n^a$ forces
a component of type $[2^m]$ for $\mfS_n^g$.}
$$   V^a_{[m,m]} \otimes V^g_{[2^m]}   \, .$$

The action of $\nabla^{ag} = \sum \pm e_{t,t}\, e_{t^\sim,t^\sim}$
restricts to the tableaux $t$ of shape $[m,m]$, and, consequently, the 
Pfaffian is proportional to 
$$   \sum_{t\in \Tab([m,m])} (-)^{\ell(\aleph,t)}\,  
       \Omega_{a,g}\,  e_{\aleph,t}^a\, e_{\aleph^\sim,t^\sim}^g \, .$$

Eventually, taking into account normalizations, one can expand the 
Pfaffian in the Young basis~:
\begin{equation}
\Pfaff\left((a_i-a_j)g[i,j]\right)  =
\sum_{t\in \Tab([m,m]} (-1)^{\ell(\zeta,t)}\,  \Y^a(t)\, \Y^g(t^\sim)  \, .
\end{equation}

In short, the Pfaffian may be considered as a trace in 
the space $ V^a_{[m,m]} \otimes V^g_{[2^m]}$.  
We summarize the preceding considerations in the following
theorem.

\begin{theorem}  \label{th:Pfaffzg}
Let $n=2m$ be an even positive integer. Let
$z[i,j]$ be Pl\"ucker indeterminates, let 
 $g[i,j]$ be indeterminates symmetrical in $i,j$, for $i,j=1\ldots n$.
Then 
\begin{eqnarray}
   \Pfaff\left( z[i,j]g[i,j]\right)  &=&
  d(\aleph)\, \Y^z(\aleph)\, \Y^g(\aleph^\sim)\, \nabla^{z,g}  \\
   &=&  \sum_{t\in \Tab([m,m])} (-1)^{\ell(\zeta,t)}\, \Y^z(t)
                \Y^g(t^\sim) \, ,
\end{eqnarray}
where $d(\aleph)$ is the number of Young tableaux of shape $[m,m]$.
\end{theorem}

This gives two ways of computing a Pfaffian. Either by a summation 
over all tableaux, or by antisymmetrization of the 
element $\Y^z(\aleph)\, \Y^g(\aleph^\sim)$ 
(one could take any other tableau than $\aleph$, or take Specht polynomials
instead of Young polynomials).

For example, for $n=6$, $z[i,j]=a_i-a_j$, $g[i,j]=(x_i^4-x_j^4)(x_i-x_j)^{-1}$,
one sees that 
$$ g[123\, |\, 456]= \Delta^x(1,2,3)\Delta^x(4,5,6) 
   S_{1,1,1}(x_1,\ldots, x_6)  \, .$$
In particular, $g[123\, |\, 456]$ is the product of a Specht polynomial
by a symmetric function in $x_1,\ldots,x_6$. 
Therefore, the Pfaffian is equal,
up to a numerical factor, to 
$$  (a_1\moins a_4)(a_2\moins a_5)(a_3\moins a_6)\Delta^x(1,2,3)\Delta^x(4,5,6)
\nabla^{ax}\, S_{1,1,1}(x_1,\ldots, x_6)  \, .$$

\section{Three symmetric groups}   

One can use $k$ families of Pl\"ucker indeterminates 
$z^1[i,j], z^2[i,j],\ldots, z^k[i,j]$,
together with a last family of  symmetric, 
or antisymmetric, indeterminates $g[i,j]$, according to the
parity of $k$. A Pfaffian 
$$\Pfaff\left( z^1[i,j]\cdots z^k[i,j]\, g[i,j]  \right)$$
   still belongs to the irreducible representation 
$$ V_{[m,m]}^{z^1} \otimes \cdots \otimes V_{[m,m]}^{z^k}    $$
of $\mfS_n^{z^1} \times  \cdots \times \mfS_n^{z^k}$.   

Therefore, it can be expanded into the Young basis~:
$$\Pfaff\left( z^1[i,j]\cdots z^k[i,j]\, g[i,j]  \right) =
\sum_{t_1,\ldots,t_k}  \Y^{z^1}(t_1)\cdots \Y^{z^k}(t_k)\,
  f(t_1,\ldots,t_k; g)   \, ,$$
sum over $k$-tuples of standard tableaux of shape $[m,m]$. 
To find the coefficients $ f(t_1,\ldots,t_k; g) $, 
one may take indeterminates $a_i^r:\, i=1\ldots n,\, r=1\ldots k$,
and put $z^1[i,j]=a_i^1-a_j^1,\ldots, z^k[i,j]=a_i^k-a_j^k$.

The main difference with the case $k=1$ is that a single specialization  
of the $a_i^r$ is not enough to determine the Pfaffian.  However,
since any element of $V^a_{[m,m]}$ is characterized by the set of all
its specializations
$$ a_{\sigma_1}=1=\cdots=a_{\sigma_m} \, ;\,
    a_{\sigma_{m+1}}=0=\cdots=a_{\sigma_n} \, , \sigma\in\mfS_n \, . $$
It is in fact sufficient to take the specializations 
corresponding to the permutations
obtained by reading the standard tableaux as permutations
(reading rows from bottom to top).
In that way, the number of specializations is equal to the number
of indeterminate coefficients, and representation theory tells us
that this system is solvable.  

I do not see anything more to say for general $k$, but shall restrict 
to $k=2$.

\begin{theorem}   \label{th:Pfaffabz}
Let $a[i,j],b[i,j],z[i,j]$, $1\leq i<j\leq n=2m$ be three families
of Pl\"ucker indeterminates. Then 
\begin{equation}
\Pfaff\left(\frac{a[i,j] b[i,j]}{z[i,j]} \right)  =
\bigl( \Y^a(\aleph) \Y^z(\zeta)\, \nabla^{az} \bigr)
\bigl( \Y^b(\aleph) \Y^z(\zeta)\, \nabla^{bz} \bigr)
\frac{d(\aleph)^2}{\prod z[i,j]}  \, ,
\end{equation} 
\end{theorem}
where $\aleph$ is the bottom tableau of shape $[m,m]$, and 
$\zeta$, the top tableau of shape $2^m$, $d(\aleph)$ still being the 
number of tableaux of the same shape as $\aleph$.

\proof As we already used, we take $a[i,j]=a_i-a_j$, $b[i,j]=b_i-b_j$.
For $z[i,j]$, we can assume that we are given a generic 
$N\times (N\plus n-1)$ matrix, $N$ sufficiently big. 
$$  M= \begin{bmatrix}    
    x_1 & x_2 &\cdots &x_n & \cdots  \\
    y_1 & y_2 &\cdots &y_n & \cdots  \\
     \cdots   & \cdots   & \cdots   & \cdots   & \fbox{$ K$}
\end{bmatrix}  \, , $$
with $K$ a submatrix of order $N\moins 2$.        
One then takes $z[i,j]$ to be the maximal minor  containing $x_i,y_j$ 
and $X_i$ (resp $Y_j$) to be the minor of order $N\moins 1$ 
containing $K$ and $x_i$ (resp. $y_j$). 
Sylvester's relation \cite{BL} states that 
$\det(K) \, z[i,j]= X_iY_j -X_j Y_i$. 
In all, the value of the Pfaffian is determined by the case
$$ \Pfaff\left( \frac{(a_i-a_j)(b_i-b_j)}{X_iY_j -X_j Y_i  }    \right)\, . $$

The first step is still to specialize
$a_1=1=\cdots=a_m$, $a_{m+1}=0=\cdots= a_n$,
 as for Th.\ref{th:Pfaffzg}.  The Pfaffian becomes
$$ F:=\det\left( (b_i-b_j)(X_iY_j -X_j Y_i)^{-1} 
             \right)_{1\leq i\leq m <j\leq n} \, .$$
One has now to specialize $b_1,\ldots,b_n$.  Instead of taking 
all permutations $\sigma\in\mfS_n$, by reordering rows and columns, 
one can suppose that there exist $\a,\b,\g$:
$\sigma= [1,\ldots,\a, \a\plus \b,\ldots, \a\plus \g, 
  \a\plus 1,\ldots, \a\plus b\moins 1, \a\plus \g\plus 1,\ldots, n]$.
In that case, the specialization 
$b_{\sigma_1}=1, \ldots, b_{\sigma_m}=1$,
 $b_{\sigma_{m+1}}=0, \ldots, b_{\sigma_n}=0$
of $F$
factorizes into two blocks, each of which
is of Cauchy type $\det( (X_iY_j-Y_iX_j)^{-1} )_{i,j=1\ldots k}$.  
For example, for $m=4$, the specialization 
$b_1=b_2=b_5=b_6=1$, $b_3=b_4=b_7=b_8=0$ 
is
$$\begin{vmatrix}
0   &   0   &   (X_1Y_7-X_7Y_1)^{-1}   &   (X_1Y_8-X_8Y_1)^{-1}   \\
0   &   0   &   (X_2Y_7-X_7Y_2)^{-1}   &   (X_2Y_8-X_8Y_2)^{-1}   \\
-(X_3Y_5-X_5Y_3)^{-1}   &   -(X_3Y_6-X_6Y_3)^{-1}   &   0   &   0   \\
-(X_4Y_5-X_5Y_4)^{-1}   &   -(X_4Y_6-X_6Y_4)^{-1}   &   0   &   0   \\
\end{vmatrix}\, .
$$
The evaluation of a Cauchy determinant is, of course, immediate (since 1812),
and in final, for any $\sigma$, the specialization 
$b_{\sigma_1}=1, \ldots, b_{\sigma_m}=1$,
 $b_{\sigma_{m+1}}=0, \ldots, b_{\sigma_n}=0$ of $F$ is
equal to
$$ (\moins 1)^{\ell(\sigma)} \Delta^z(\sigma_1,\ldots,\sigma_m)
  \Delta^z(\sigma_{m+1},\ldots,\sigma_n)\,
 \prod_{1\leq i\leq m <j\leq n} (z_i-z_j)^{-1} \, .$$
Therefore, $F$ coincides, up to a numerical factor, with        
\begin{multline*}
b_1\cdots b_m \Delta^z(1,\ldots,m)\Delta^z(m\plus1,\ldots,n)  \nabla^{b,z}
 \prod_{1\leq i\leq m <j\leq n} (z_i-z_j)^{-1}  \\
= (b_1-b_{m+1})\cdots (b_m-b_n) \Delta^z(m\plus1,\ldots,n)  \nabla^{b,z}
 \prod_{1\leq i\leq m <j\leq n} (z_i-z_j)^{-1}  \, .
\end{multline*}
From this specialization, one writes the Pfaffian as
$$ \Y^b(\aleph) \Y^z(\zeta)\, \nabla^{bz}
  \frac{\Delta^z(1,\ldots,m) \Delta^z(m\plus 1,\ldots,n)}
 {\Delta^z(1,\ldots,n)} \,  \nabla^{abz} \, . $$
Since any permutation in $\mfS^{abz}$ commutes with
 $\nabla^{bz}\Delta^z(1,\ldots,n)^{-1}$,
the theorem follows  \QED

Okada \cite[Th.3.4]{Okada} (see also    
\cite[Formula 1.8]{IshiOkada})
has already computed 
$\Pfaff\bigl( (a_i\moins a_j)(b_i\moins b_j) (x_i\moins x_j)^{-1} \bigr)$. 
His formula can be written    
\begin{multline}  \label{Pfaffabx}
  \Pfaff\left(\frac{(a_i-a_j)(b_i-b_j)}{x_i^2-x_j^2}\right) \\ =
 \prod_{i<j}\frac{x_i+x_j}{x_i-x_j}\,
\Pfaff\left(\frac{a_i-a_j}{x_i+x_j} \right)
 \, \Pfaff\left(\frac{b_i-b_j}{x_i+x_j} \right)  \, .
\end{multline}

In \cite{Okada} and \cite{IshiOkada}, one finds many evaluations
of Pfaffians and determinants, with entries which are specializations
of Pl\"ucker indeterminates.   For example, Okada \cite[Th.3.4]{Okada}
takes 
\begin{equation*}
a[i,j]= \begin{vmatrix} 1+ a_i x_i & x_i +a_i \\  
            1+ a_j x_j & x_j +a_j\end{vmatrix}
\ \text{or}\ 
a[i,j] =
\begin{vmatrix}
 1+ a_i x_i^2 & x_i +a_i x_i  & x_i^2 +a_i\\
 1+ a_j x_j^2 & x_j +a_j x_j  & x_j^2 +a_j\\
 1+cz^2       & z+cz          & z^2 +c
\end{vmatrix}  \, .
\end{equation*}

All these families satisfy the Pl\"ucker relations, and one could take
determinants of higher order of the type below, as indeterminates $z[i,j]$.
On the other hand, Pfaffians involving elliptic functions as
in \cite{Okada2} do not fall in this category, 
the Riemann relations replacing  in that case
Pl\"ucker relations.

\section{Special Pfaffians}    

There are cases where  
$\Y^z(\aleph)\, \Y^g(\aleph^\sim)\, \nabla^{z,g}$ 
can be written as a determinant.
Indeed, let us take $z[i,j]=(a_i-a_j)$, and $g[i,j]=(x_i+x_j)^{-1}$
as Sundquist\cite{Sundquist}. 
For any integer $k$, let $U(a,x^k)$ be the determinant of oder $n$
\begin{equation}  \label{UU}
  U(a,x^k) :=\Bigl|a_i x_i^0,\, a_i x_i^k,\dots,
  a_i x_i^{k(m-1)},\, x_i^0,\, x_i^k,\ldots, x_i^{k(m-1)}
       \Bigr|_{i=1\ldots n}\ .
\end{equation}

The Laplace expansion of $U(a,x)$  
along its first $m$ columns shows that it belongs to the 
Specht representation of
$\mfS_n^x$ of shape $[2^m]$, and therefore, that $U(a,x)$ is an 
element of $V_{[m,m]}^a\otimes V_{[2^m]}^x$, 
which is identified by the specialization 
$a_1=1 = \cdots = a_m$,
 $a_{m+1} =0=\cdots = a_n$.
Needless to add that this specialization is 
$$\Delta(x_1,\ldots,x_m)\, \Delta(x_{m+1},\ldots, x_n) \, .$$

On the other hand, the Pfaffian  
$\Pfaff\left( (a_i\moins a_j)(x_i\plus x_j)^{-1}  \right)$ 
belongs to the same space, and the same 
specialization sends it, according to (\ref{Detg}),
and thanks to the Cauchy identity relative to the determinant  
$\Bigl|(x_i+y_j)^{-1} \Bigr|_{i,j=1\ldots m}$, to
$$ c_{\aleph}(x) = \frac{\Delta^x(1,\ldots,m)\,
  \Delta^x(m\plus 1,\ldots, n)  }
       {\prod_{i\leq m <j} x_i+x_j }\, .  $$

Taking  a   symmetrical denominator, one rather writes
$$ c_{\aleph}(x) = \frac{\Delta^{xx}(1,\ldots,m)\,
  \Delta^{xx}(m\plus 1,\ldots, n)  }
       {\prod_{1\leq i<j\leq n} x_i+x_j }\, .  $$

Therefore, this specialization coincides with the one of $U(a,x^2)$,
and one recovers the following theorem.

\begin{theorem}[Sundquist]    \label{th:Sundquist}
Let $n=2m$, $a_1,\ldots,a_n$, $x_1,\ldots,x_n$ be indeterminates.
Let $U(a,x)$ be the determinant \em\eqref{UU}.

Then   
\begin{eqnarray*}
   \Pfaff\left( \frac{a_i-a_j}{x_i+x_j} \right) 
  &=& \frac{d}{n!} \sum_\sigma  (\moins 1)^{\ell(\sigma)}
       \Bigl( (a_1\moins a_{m+1})\cdots (a_{m-1}\moins a_n)\, \times \\
 & & \hspace{30pt}
\Delta^{xx}(1,\ldots,m) \Delta^{xx}(m\plus 1,\ldots, n) )\Bigr)^\sigma
   \, \prod_{1\leq i<j\leq n} (x_i+x_j)^{-1}   \\
 &=&  \det\bigl(U(a,x^2) \bigr)\, 
\prod_{1\leq i<j\leq n} (x_i+x_j)^{-1} \, ,
\end{eqnarray*}
where the sum is over all permutations $\sigma\in\mfS_n^{a,x}$,
and $d$ is the number of Young tableaux of shape $[m,m]$.
\end{theorem}

We have given another expression in Th\ref{th:Pfaffzg}, using the Young basis.
To stay nearer the expansion of $U(a,x)$, we  
 can also use the Specht basis,
but there will be extra terms corresponding to pairs of non-orthogonal
tableaux (as soon as $n=6$).
 
Indeed, for $n=4$, we have
\begin{multline*}
\prod_{i<j}(x_i+x_j)\, 
 \Pfaff\left( \frac{a_i-a_j}{x_i+x_j} \right) = 
        (a_1-a_2)(a_3-a_4)(x_1^2 -x_3^2)(x_2^2-x_4^2)  \\
   -(a_1-a_3)(a_2-a_4)(x_1^2 -x_2^2)(x_3^2-x_4^2)   \, .
\end{multline*}

For $n=6$, the expansion in the Specht basis is 
\begin{multline*}
\prod_{i<j}(x_i+x_j)\, 
\Pfaff\left( \frac{a_i-a_j}{x_i+x_j} \right) =
  (a_1-a_2)(a_3-a_4)(a_5-a_6) \times \\
   \Delta(x_1^2,x_3^2,x_5^2) \Delta(x_2^2, x_4^2, x_6^2)  
  \bigl(1 -s_2^{ax} -s_4^{ax} +s_2^{ax}s_4^{ax} -
   s_2^{ax}s_4^{ax}s_3^{ax} \bigr)   \\
 - (a_1-a_2)(a_3-a_4)(a_5-a_6) \, \Delta(x_1^2,x_2^2,x_3^2)
    \Delta(x_4^2, x_5^2, x_6^2) \\
=\begin{vmatrix}
\, a_1 & a_1 x_1^2 &a_1 x_1^4 & 1 & x_1^2 & x_1^4 \, \\
\, a_2 & a_2 x_2^2 &a_2 x_2^4 & 1 & x_2^2 & x_2^4 \, \\
\, a_3 & a_3 x_3^2 &a_3 x_3^4 & 1 & x_3^2 & x_3^4 \, \\
\, a_4 & a_4 x_4^2 &a_4 x_4^4 & 1 & x_4^2 & x_4^4 \, \\
\, a_5 & a_5 x_5^2 &a_5 x_5^4 & 1 & x_5^2 & x_5^4 \, \\
\, a_6 & a_6 x_6^2 &a_6 x_6^4 & 1 & x_6^2 & x_6^4 
\end{vmatrix} \, .
\end{multline*} 
The first five terms, written as images of the first one,
 are   the Specht polynomials for pairs of tableaux
transposed of each other, but there is a sixth term corresponding
to the only  non-zero entry outside the diagonal in the matrix
of scalar products.

In the case where the $a_i's$ are fixed powers of the indeterminates $x_i$'s,
then the determinant $U(a,x^2)$ is a determinant of powers of $x_i$'s,
proportional to a Schur function.
Thus, Th.\ref{th:Sundquist} implies
\begin{corollary}  \label{cor:PfafPowers}
Let $r,k$ be positive integers, $n=2m$ be an even integer. 
Let $q=2(k-1)$, and $\lambda$ be the (increasing) partition
$[0,r,q, q\plus r, \ldots, (m\moins 2)q, (m\moins 2)q+r ,
(m\moins 1)q, m\moins 1)q+r]$
Then 
$$  \Pfaff\left( \frac{x_i^{r+1}-x_j^{r+1}}{x_i^k+x_j^k} \right)   =
\frac{\Delta^x(1,\ldots,n)}{\prod_{1\leq i<j\leq n} x_i^k+x_j^k}\ 
 S_\l(x_1,\ldots, x_n)   \, .  $$    
\end{corollary}   

We can give more interesting examples of Pfaffians 
$\Pfaff\left( z[i,j]g[i,j]\right)$, with $g[i,j]$ a symmetric function
in $x_i,x_j$, which admits a symmetric function in $x_1,\ldots, x_n$
as a factor.

For example, take a partition $\lambda$, an alphabet $B$,
and variables $x_1,\ldots,x_n$.

Chosing a positive  $k$, we want to evaluate 
$$ \Pfaff\left((a_i-a_j) S_\l(B+ x_i\plus x_j)\frac{x_i-x_j}{x_i^k-x_j^k}
   \right)  \, . $$
According to Th.\ref{th:Pfaffzg}, we need only compute the
specialization $a_1=1,\ldots,a_m=1$, $a_{m+1}=0,\ldots, a_n=0$,
which is equal to 
$$  \det\left({S_\lambda(B+x_i\plus x_j)}\frac{x_i-x_j}{x_i^k -x_j^k}  
         \right)_{ 1\leq i\leq m <j\leq n} \, .$$   

To proceed further, one needs to evaluate such determinants.
We shall do that in the next section. For the moment, let us only  
use the fact that the determinant in question is \emph{the product of a 
symmetric function} $f(x_1,\ldots, x_n)$ by 
$\Delta^x(1\ldots m)\Delta^x(m\plus 1\ldots n) 
 \prod_{1\leq i\leq m <j\leq n} (x_i\moins x_j) (x_i^k \moins x_j^k)^{-1} $.

The symmetric function can then  be factored out, so that 
the Pfaffian 
$\Pfaff\left((a_i-a_j) S_\l(B+ x_i\plus x_j)
   (x_i\moins x_j) (x_i^k \moins x_j^k)^{-1} \right)$
is finally equal to
\begin{equation}    \label{PfaffSf}  
  f(x_1,\ldots, x_n)\, \Pfaff\left( (a_i-a_j) 
                       \frac{x_i-x_j}{x_i^k -x_j^k}\right) \, . 
\end{equation}
Thus, the evaluation of the Pfaffian of order $2m$
has been reduced to the evaluation of a determinant 
of order $m$. 

Ishikawa \cite{Ishikawa}, Okada \cite{Okada},
  and Ishikawa-Okada-Tagawa-Zeng \cite{IshiOkada}
have given many generalizations of Sundquist's Pfaffian.
Instead of using Pl\"ucker coordinates, they use specific determinants
(which, of course, satisfy built-in Pl\"ucker relations).

\section{Determinants and two symmetric groups}

The fundamental surveys of Krattenthaler \cite{Kratt,Kratt2} 
describe many methods
to evaluate determinants. We would like to add to them one more method,
using two symmetric groups. 

In the course of proving Th.\ref{th:Pfaffabz}, we have 
met a determinant which 
 happened to possess an unsuspected  global symmetry, and that we record
now (notice that $\nabla^{az}$ is given by a summation on the full
symmetric group $\mfS_n$).

\begin{corollary}   \label{cor:Detaz}
Let $\{ a[i,j]\}$, $\{ z[i,j]\}$, $1\leq i <j\leq n=2m$
 be two families of Pl\"ucker indeterminates. Then 
\begin{equation}
\det\left(  \frac{ a[i,j]}{z[i,j]}
\right)_{1=1\ldots m,\, j=m\!+\!1\ldots n}
= d(\aleph)\, \Y^a(\aleph) \Y^z(\zeta) \, \nabla^{az} \, \prod_{1\leq i <j\leq n}
   z[i,j]^{-1}   \, ,
\end{equation}
where $\aleph$ is the bottom tableau of shape $[m,m]$, $\zeta$,
the top tableau of shape $2^m$.
\end{corollary}

The special case where $z[i,j]=x_i-x_j$, or $z[i,j]=x_i^k-x_j^k$
is worth commenting, since it reveals a symmetry in $x_1,\ldots, x_n$
that we emphasize in the next theorem theorem.

To evaluate 
$\det\bigl( a[i,j] (x_i-x_j)^{-1}\bigr)_{1=1\ldots m,\, j=m\!+\!1\ldots n}$,
 we already used that we need only take $a[i,j]=a_i-a_j$, and specialize 
$[a_1,\ldots, a_n]$ in all permutations of $[1^m,0^m]$. 
However, writing 
$$ R^x(1\ldots m\, |\, m\plus 1\ldots n) := \prod\nolimits_{1\leq i\leq m}
   \prod\nolimits_{m+1\leq j\leq n}   (x_i-x_j)  \, ,$$
 it is clear that 
$$ R^x(1\ldots m\, |\, m\plus 1\ldots n)\,
       \Biggl|\frac{a_i-a_j}{x_i-x_j} \Biggr|_{1\leq i\leq m<j\leq 2m} $$
has the same specializations as  $U(a,x)$. 

Thanks to Th.\ref{th:Pfaffzg} and Th.\ref{th:Sundquist},
going back to the variables $z[i,j]$, taking variables $x_i^2$ instead
of $x_i$, we have just obtained~:

\begin{theorem}  \label{th:DetSf}
 Let $z[i,j]$, $1\leq i<j \leq n=2m$ be Pl\"ucker indeterminates,
and let $x_1,\ldots, x_n$ be indeterminates.  Then 
\begin{eqnarray}
  \notag
  \det\left| \frac{z[i,j]}{x_i^2-x_j^2}    \right|_{1\leq i\leq m<j\leq n}
  &=&   \frac{1}{R^{xx}(1\ldots m\, |\, m+1\ldots n) }   \\
 & & \hspace{15pt} \times \
  \sum_{t\in \Tab([m,m]} (-1)^{\ell(\zeta,t)}\, \Y^z(t) \Y^{xx}(t^\sim)  \\ 
  &=& \Y^z(\aleph) \Y^{xx}(\aleph^\sim)\, \nabla^{xz}
            \frac{d(\aleph)}{R^{xx}(1\ldots m | m\plus1\ldots n) }   \\
  &=& \frac{\prod_{1\leq i<j\leq n}{x_i+x_j}}{
                R^{xx}(1\ldots m\, |\, m+1\ldots n) } \,   
         \Pfaff\left(\frac{z[i,j]}{x_i+x_j}  \right)  \, .    
\end{eqnarray}
the Pfaffian being of order $n$, and the superscript $xx$ meaning using the
indeterminates $x_i^2$ instead of $x_i$. 

In particular, $R^x(1\ldots m|m\plus 1\ldots n)\, 
                  \left| \frac{z[i,j]}{x_i-x_j} \right|$
is symmetrical in $x_1,\ldots, x_n$.
\end{theorem}

Taking $z[i,j]= S_\lambda(B+x_i+x_j)(x_i-x_j)$
and changing the powers of the variables $x_i$ in denominator, 
we get the following corollary.

\begin{corollary}   \label{th:DetSf2}
Let $\lambda$ be a partition, $B$ be an alphabet, $m,k$ be two positive
integers.  Then
$$ \frac{1}{\Delta^x(1\ldots m)\Delta^x(m\plus 1\ldots 2m)}
 \left(\prod_{1\leq i\leq m<j \leq 2m} \frac{x_i^k-x_j^k}{x_i-x_j}\right)\,
  \left| \frac{S_\lambda(B \plus x_i\plus x_j)}{
    S_{k-1}(x_i+x_j)}    \right|_{1\leq i\leq m<j\leq 2m}
$$
is a function symmetrical in $x_1,\ldots, x_{2m}$. 
\end{corollary}

Notice that $\Bigl|p_2(B+x_i+x_j)  \Bigr|_{1\leq i\leq 2<j\leq 4}$,
where $p_2$ is the second power sum, does not furnish a symmetric function
in $x_1,\ldots, x_4$ (this does not contradict the corollary, because,
fortunately, $p_2$ is not a Schur function).

For $k=1$, and $\lambda=r^p$, a rectangular partition, one gets 
an identity which is useful in the theory of orthogonal
polynomials \cite[Prop. 8.4.3]{Cbms}~:
\begin{multline}
\frac{1}{ \Delta^x(1\ldots m)\Delta^x(m\plus 1\ldots 2m)} 
  \left|  S_{r^p}(B+x_i+x_j)\right|_{1\leq i\leq m<j\leq n} \\
=  \left(S_{(r+1)^{p-1}}(B)  \right)^{m-1}\, 
   S_{(r-m+1)^{p+m-1}}(B+x_1+\cdots +x_n)  \, .
\end{multline}

More precisely, given \emph{moments} $\mu_k=(\moins 1)^k S_{1^k}(B)$,
supposed to be sufficiently generic, 
then $\{S_{n^n}(B+x)  \}$ is a family of orthogonal polynomials, with
respect to this moments \cite[Ch.8]{Cbms}.   
The Christoffel-Darboux kernel of order $r$,
$K_r(x,y)$, is proportional to $S_{r^{r+1}}(B+x+y)$,
and Cor.\ref{th:DetSf2}  states that the determinant 
with entries $K_r(x_i,y_j)$, $1\leq i,j\leq n$ is a symmetric function
of $x_1,\ldots,x_n,y_1,\ldots,y_n$ and gives its precise value.
This property can be directly proved, using Bazin relation on minors
\cite{L-Shi} (see also the article of Rosengren 
about the relations between Pfaffians and kernels \cite{Rosengren}).

The case  $k=2$, and $\lambda=\rho :=[1,\ldots,r]$, 
has been settled by \cite{IshiOkada, Okada2}, and reads
\begin{multline}  
\left| \frac{S_\rho(B+x_i+x_j)}{x_i+x_j}\right|_{1\leq i\leq m<j\leq n}
  \\ =
  \frac{ \Delta^x(1\ldots m)\Delta^x(m\plus 1\ldots 2m) }{
     \prod_{1\leq i\leq m<j \leq 2m}  x_i+x_j}  
   \left(S_\rho(B)  \right)^{m-1}\,
   S_\rho(B+x_1+\cdots +x_n)  \, .
\end{multline}

In these two cases, the symmetric function has been further factorized,
compared to the case of a general partition $\lambda$.
This induces a factorization of Pfaffians~:
\begin{multline}
 \Pfaff\left(a[i,j]S_{r^p}(B+x_i+x_j) \right)   \\
  =  \left(S_{(r+1)^{p-1}}(B)  \right)^{m-1}\,
         S_{(r-m+1)^{p+m-1}}(B+x_1+\cdots +x_n) \hspace*{45pt}  \\ 
  \Y^a(\aleph) \Delta^x(1\ldots m)
          \Delta^x(m\plus 1\ldots 2m)    \nabla^{ax}   \, ,
\end{multline}
\begin{multline}
   \Pfaff\left( \frac{ a[i,j]}{x_i+x_j}\, S_\rho(B+x_i+x_j) \right) \\
  = \left(S_\rho(B)  \right)^{m-1}\,
   S_\rho(B+x_1+\cdots +x_n)  \,  
           \Pfaff\left( \frac{a[i,j]}{x_i+x_j}\right)	\, .
\end{multline}

Our last example will be related to the six-vertex model in physics.
Stroganov \cite{Stroganov} found that the determinant with entries
$$   \sin( x_i-y_j +\eta)^{-1} \sin( x_i-y_j -\eta)^{-1}  \, ,$$
with $\eta=exp(\pi\sqrt{-1}/3)$, $i,j=1\ldots n$,
is the product of a symmetric function in $x_1,\ldots,x_n,y_1,\ldots,y_n$
by $\Delta^x(1,\ldots,n) \Delta^y(1,\ldots,n)$.   

Since 
$$ -4\sin( x_i-y_j +\eta) \sin( x_i-y_j -\eta)=
  1+ 2 \cos(2(x\moins y))  = 1 +a/b + b/a , $$
with $a= exp(2\pi x\sqrt{-1})$, $b=exp(2\pi y\sqrt{-1})$,
Stroganov's case is the evaluation of the determinant
$\Bigl| \bigl( x_i/y_j-y_j/x_i  \bigr)\bigl(
   (x_i/y_j)^3 -(y_j/x_i)^3   \bigr)^{-1}   \Bigr|$.

The following lemma gives a more general case, as a corollary of 
Th. \ref{th:DetSf}.

\begin{lemma}
Let $k,r$ be two positive integers. Then 
$$ \det\Bigl( (x_i^r-x_j^r)(x_i^k-x_j^k)^{-1}
                  \Bigr)_{i=1\ldots m,\, j=m+1\ldots 2m}$$    
is equal to the product of the Schur function in $x_1,\ldots, x_{2m}$
of index
$$ [0,\g,\b, \b\plus \g,2\b, 2\b\plus \g,\ldots, (n\moins 1)\b, 
(n\moins 1)\b\plus \g]
$$
times $\Delta^x(1\ldots m) \Delta^x(m\plus 1\ldots 2m)$, where
$\g=r-1$, $\b=k-2$.
\end{lemma}  

For example, for $m=3$, $r=2$, $k=5$,  then
the determinant is equal to
$$\Delta^x(1,2,3)\Delta^x(4,5,6)\,  s_{[0,1,3,4,6,7]}(x_1,\ldots, x_6)\, .$$

\section{Note: Young's basis}

Young first defined \emph{natural idempotents}, giving rise to what 
we have called the \emph{Specht basis}. He then obtained orthogonal idempotents by an orthogonalization process which was later clarified by
Thrall (see Rutherford \cite{Rutherford}). 

The easiest way of obtaining Young's orthogonal idempotents 
is to characterize them as simultaneous eigenvectors for the 
\emph{Jucys-Murphy elements} 
$$  \xi_j := \sum_{i<j} (i,j) \ , \ j=0,\ldots,n \, ,$$
where the sum is over transpositions (cf. Okounkov-Vershik \cite{OV}).

However, this approach does not provide the relations between the different idempotents for the same shape, and is inappropriate for our decomposition
of Pfaffians.

We need to reinterpret Young's orthogonalization in terms of
the \emph{Yang-Baxter relations}~:
\begin{gather}  \label{eq:Yang1}
  \left( s_i+\frac{1}{\a}\right) 
  \left(s_{i+1}+\frac{1}{\a \plus\b}\right) \left(s_i+\frac{1}{\b}\right) =
\left( s_{i+1}+\frac{1}{\b}\right) \left(s_{i}+\frac{1}{\a\plus \b}\right) 
\left(s_{i+1}+\frac{1}{\a}\right)  \\ \label{eq:Yang2}
  \left( s_i+\frac{1}{\a}\right)\, \left(s_j+ \frac{1}{\b}\right)  = 
 \left(s_j+ \frac{1}{\b}\right)\,
  \left( s_i+\frac{1}{\a}\right)  \ ,\ |i-j|\neq 1
\end{gather}
The graphical representation of these relations
is easy to remember (taking $i=1$)~:
$$ \hspace{-1cm} \begin{array}{rcl}
   & [1\,2\,3]    &\\[2pt]
   \hbox{${\hskip -4mm\raise 3mm
    \hbox{$ s_1+\frac{1}{\a} $}\, \diagup\hskip -2mm\diagup }$}  &
   & \diagdown \raise 3mm \hbox{$\, s_2+\frac{1}{\b}  $}  \\[5pt]
   [2\,1\,3]  &  & [1\,3\,2]  \\[5pt]
   \hbox{$ s_2+\frac{1}{\a +\b}  \, \bigg|$\hskip 3.7mm}  &
    &\hbox{\hskip 3.4mm $\bigg\| \, s_1+\frac{1}{\a+\b} $}  \\[2pt]
   [2\,3\,1]  & &[3\,1\,2] \\[5pt]
   \hbox{$ s_1+\frac{1}{\b} \ \diagdown\hskip -2mm\diagdown$} &
      &\diagup \ s_2+\frac{1}{\a} \\[5pt]
  &[3\,2\,1] &
\end{array} \hskip 15mm
 \begin{array}{rcl}
   & [1\,2\,3\,4]    &\\[2pt]
   \hbox{$\hskip -4mm\raise 3mm  \hbox{$ s_1+\frac{1}{\a} $}\, \diagup\hskip
                                  -2mm\diagup $}     &
  & \diagdown \raise 3mm\hbox{$\, s_3+\frac{1}{\b}$}  \\[5pt]
   [2\,1\,3\,4] & &[1\,2\,4\,3]  \\[5pt]
   s_3+\frac{1}{\b} \ \diagdown & & \diagup\hskip -2mm\diagup \  s_1+\frac{1}{\a} \\[5pt]
  &[2\,1\,4\,3] &
\end{array} $$

The standard Young tableaux of a given shape are the vertices of
a graph obtained by generating them with simple transpositions,
starting from the top  $\zeta$. We keep the same directed
graph, but label an edge \hbox{$t \to t s_i$}  by
$s_i+1/\rho$~:
$$  t \ \stackrel{s_i+1/\rho}{\hbox to 33pt{\rightarrowfill}}
                \   t\, s_i \, , $$
where $\rho$ is the diagonal distance (the difference of \emph{contents})
between the letters $i$ and $i\plus 1$ in $t$.  

A path in such a graph is interpreted as the product, in the group algebra,
of the edges composing it, and the Yang-Baxter relations insure that 
two paths having the same end points evaluate to the same element
in the group algebra.

We replace now $\zeta$ by the Specht polynomial $\Delta_\zeta^a$,
and define, for any other standard tableau of the same shape,
the \emph{Young polynomial} $\Y^a(t)$ by~: 
$$  \Y^a(t) = \Delta_\zeta^x\, \left(s_i +\frac{1}{\rho}\right)
       \cdots \left( s_j+\frac{1}{\rho'}\right)  \, ,$$
whenever 
$$  \zeta \ \stackrel{s_i+1/\rho}{\hbox to 33pt{\rightarrowfill}}  
  \ \cdots \stackrel{s_j+1/\rho'}{\hbox to 33pt{\rightarrowfill}}  
\ t $$
is a path from $\zeta$ to $t$.

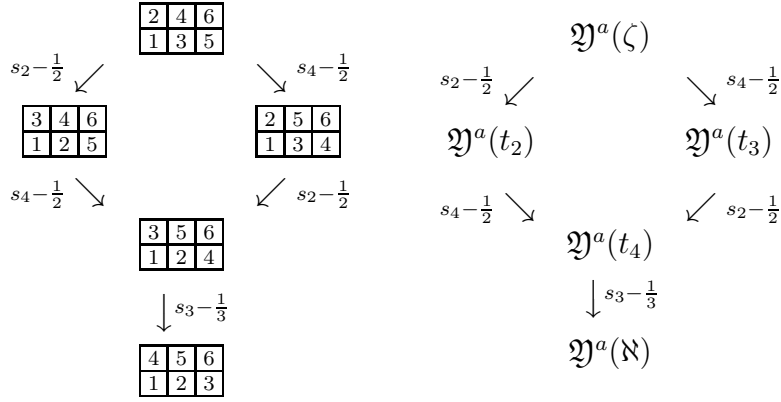
\begin{equation*}
\begin{array}{rcl}
          & \smallyoung{\s 2&\s4&\s6\cr \s1&\s3&\s5\cr}                \\[4pt]
  \raise 5pt \hbox{$\s s_2 -\frac{1}{2}$} \swarrow &
            & \searrow  \raise 5pt \hbox{$\s s_4-\frac{1}{2}$} \\{}
\smallyoung{\s3&\s4&\s6\cr \s1&\s2&\s5\cr} & & 
 \smallyoung{\s2&\s5&\s6\cr \s1&\s3&\s4\cr}\\[10pt]
{\s s_4-\frac{1}{2}} \searrow &&
   \swarrow {\s s_2 -\frac{1}{2}} \\
               & \smallyoung{\s3&\s5&\s6\cr \s1&\s2&\s4\cr} &  \\[13pt]
  &   \phantom{\s s_3} \big\downarrow 
    \raise 3pt \hbox{$\s s_3-\frac{1}{3} $} &\\[2pt]
                       & \smallyoung{\s4&\s5&\s6\cr \s1&\s2&\s3\cr}
 \end{array}\quad\quad  
\begin{array}{rcl}
          & \Y^a(\zeta)                \\[4pt]
  \raise 5pt \hbox{$\s s_2 -\frac{1}{2}$} \swarrow &
            & \searrow  \raise 5pt \hbox{$\s s_4-\frac{1}{2}$} \\[5pt]
\Y^a(t_2) & &  \Y^a(t_3)\\[10pt]
{\s s_4-\frac{1}{2}} \searrow &&
   \swarrow {\s s_2 -\frac{1}{2}} \\
               & \Y^a(t_4) &  \\[6pt]
           & \phantom{\s s_3} \big\downarrow 
      \raise 3pt \hbox{$\s s_3-\frac{1}{3} $} &\\[6pt]
                       & \Y^a(\aleph)
 \end{array}
\end{equation*}
The graph on the right side describes the generation of Young's basis,
starting with $\Y\left( \begin{smallmatrix} 
2&4&6\\ 1&3&5 \end{smallmatrix}\right)$, and applying 
those $s_i +1/\rho$ which label the edges.

Our graph is directed, but since $-\rho$ is the distance between 
$i,i\plus1$ in $ts_i$, if $\rho$ is the distance in $t$, and since
$$ (s_i+1/\rho)\, (s_i -1/\rho) = 1-1/\rho^2 \, ,  \rho \neq 1\, , $$
one could use a double orientation by normalizing the edges,
taking\\  $(s_i+1/\rho)/\sqrt{ 1-1/\rho^2}$ instead of
$(s_i+1/\rho)$.

Pfaffians are obtained by taking a space $V_{[m,m]}^a\otimes V_{[2^m]}^x$,
and using  a pair of orthonormal bases that we write below 
 in terms of the two Young bases. The graph on the left describes
the orthonormal basis for shape $[3,3]$, generated downwards,
and the graph on the right,
the basis for shape $[2,2,2]$, generated upwards. The Pfaffian is obtained
by taking the sum of products of the corresponding 
vertices of the two graphs~:
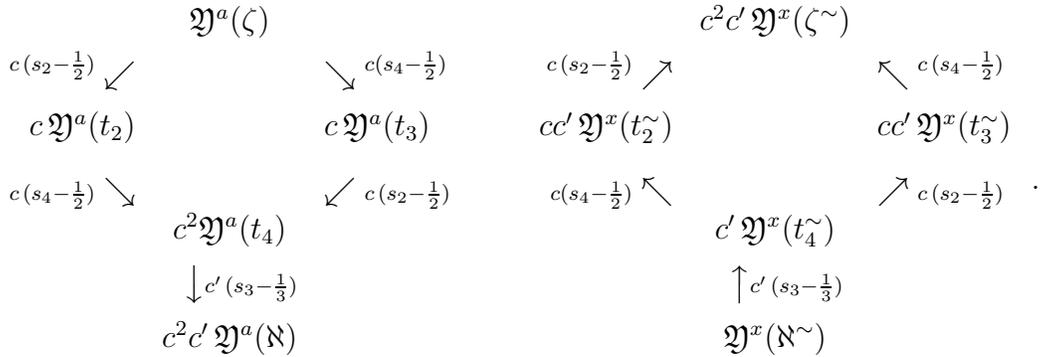
\begin{equation*}
\begin{array}{rcl}
          & \Y^a(\zeta)                \\[4pt]
  \raise 5pt \hbox{$\s c\,(s_2 -\frac{1}{2})$} \swarrow &
            & \searrow  \raise 5pt \hbox{$\s c(s_4-\frac{1}{2})$} \\[5pt]
       c\, \Y^a(t_2) & &  c\, \Y^a(t_3)\\[10pt]
{\s c\, (s_4-\frac{1}{2})} \searrow &&
   \swarrow {\s c\, (s_2 -\frac{1}{2})} \\
               & c^2\Y^a(t_4) &  \\[6pt]
           & \phantom{\s s_3} \big\downarrow{\s c'\, (s_3-\frac{1}{3})}&\\[6pt]
                       &c^2c'\, \Y^a(\aleph)
 \end{array}
\quad\quad
\begin{array}{rcl}
          & c^2c'\, \Y^x(\zeta^\sim)                \\[4pt]
  \raise 5pt \hbox{$\s c\, (s_2 -\frac{1}{2})$} \nearrow &
            & \nwarrow  \raise 5pt \hbox{$\s c\, (s_4-\frac{1}{2})$} \\[5pt]
cc'\, \Y^x(t_2^\sim) & &  cc'\, \Y^x(t_3^\sim)\\[10pt]
{\s c(s_4-\frac{1}{2})} \nwarrow &&
   \nearrow {\s c\, (s_2 -\frac{1}{2})} \\
               & c'\, \Y^x(t_4^\sim) &  \\[6pt]
           & \phantom{\s s_3} \big\uparrow {\s c'\, (s_3-\frac{1}{3})} &\\[6pt]
                       & \Y^x(\aleph^\sim)
 \end{array} \ .
\end{equation*}
The normalization  constants $\frac{1}{\sqrt{1-1/\rho^2}}$ 
are 
$c= \frac{1}{\sqrt{1-1/4}}$ and $c'=\frac{1}{\sqrt{1-1/9}}$, 
because the diagonal 
distances involved are $\pm 2$ and $\pm 3$.

In conclusion, the sum 
$\sum_{t\in\Tab([3,3])} \pm \Y^a(t) \Y^x(t^\sim)$ is equal,
up to the global factor $c^2c'$, to the 
same sum when using the orthonormal basis
rather than the Young polynomials (which are, as we chose to define them,
only an orthogonal basis).

Notice that the sum of products of Young polynomials can be written
\begin{multline*}
  \Y^a(\aleph)\Y^x(\aleph^\sim)\, 
  \Biggl( 1-(s_3^a+{\s\frac{1}{3}})(s_3^x-{\s\frac{1}{3}})
   \frac{1}{1-\frac{1}{9}}  \times \\ 
  \left( 1-(s_2^a+{\s\frac{1}{2}})(s_2^x-{\s\frac{1}{2}})
       \frac{1}{1-\frac{1}{4}}  \right)
   \left( 1-(s_4^a+{\s\frac{1}{2}})(s_4^x-{\s\frac{1}{2}})
       \frac{1}{1-\frac{1}{4}}  \right)
   \Biggr)\\
= \quad \Y^a\left(\begin{smallmatrix}  4&5&6\\ 1&2&3 \end{smallmatrix}\right)
  \Y^x\left(\begin{smallmatrix} 3&6\\ 2&5\\ 1&4\end{smallmatrix}\right)
 - \Y^a\left(\begin{smallmatrix}  3&5&6\\ 1&2&4 \end{smallmatrix}\right)
  \Y^x\left(\begin{smallmatrix} 4&6\\ 2&5\\ 1&3\end{smallmatrix}\right)
   \times \\
\left( 1-(s_2^a+{\s\frac{1}{2}})(s_2^x-{\s\frac{1}{2}})
       \frac{1}{1-\frac{1}{4}}  \right)
   \left( 1-(s_4^a+{\s\frac{1}{2}})(s_4^x-{\s\frac{1}{2}})
       \frac{1}{1-\frac{1}{4}}  \right) \\
 = \cdots = \quad
  \Y^a\left(\begin{smallmatrix}  4&5&6\\ 1&2&3 \end{smallmatrix}\right)
  \Y^x\left(\begin{smallmatrix} 3&6\\ 2&5\\ 1&4\end{smallmatrix}\right)
 - \Y^a\left(\begin{smallmatrix}  3&5&6\\ 1&2&4 \end{smallmatrix}\right)
  \Y^x\left(\begin{smallmatrix} 4&6\\ 2&5\\ 1&3\end{smallmatrix}\right) \\
 + \Y^a\left(\begin{smallmatrix}  3&4&6\\ 1&2&5 \end{smallmatrix}\right)
  \Y^x\left(\begin{smallmatrix} 5&6\\ 2&4\\ 1&3\end{smallmatrix}\right)
 +\Y^a\left(\begin{smallmatrix}  2&5&6\\ 1&3&4 \end{smallmatrix}\right)
  \Y^x\left(\begin{smallmatrix} 4&6\\ 3&5\\ 1&2\end{smallmatrix}\right)
-\Y^a\left(\begin{smallmatrix}  2&4&6\\ 1&3&5 \end{smallmatrix}\right)
  \Y^x\left(\begin{smallmatrix} 5&6\\ 3&4\\ 1&2\end{smallmatrix}\right)
\end{multline*}

In the preceding sections,  we did not have recourse to normalization constants,
but used the Young basis and checked the overall factor by
computing a specialization of the Pfaffian.

\bigskip
\noindent
\emph{Acknowledgment.  
The author benefits from the ANR project BLAN06-2\_134516. 
This work was partly done  in January-February 2005,
during a Combinatorial Semester 
at the Mittag-Leffler Institute.
I thank the Institute for its warm hospitality.}

\end{document}